\documentclass[11pt]{article}
\usepackage{amssymb,amsmath}
\usepackage[T1]{fontenc}
\usepackage[utf8]{inputenc}
\usepackage{graphicx}
\usepackage{amscd}
\usepackage{epsfig}
\usepackage{mathrsfs}
\usepackage[usenames,dvipsnames]{xcolor}
\usepackage{cite}
\usepackage{epigraph}
\begin{document}

\begin{center}
{\bf \Large A novel approach to solving a multipoint boundary value problem for an integro-differential equation}
\end{center}

\begin{center}
\textbf{ Anar T. Assanova$^1$, Elmira A. Bakirova$^2$, Roza E. Uteshova$^3$}
 \\ [0.2cm]
{\it Institute of Mathematics and Mathematical Modeling, \\
 125, Pushkin Str., 050010, Almaty, Kazakhstan }\end{center}
$~^1$e-mail: assanova@math.kz, $~^2$e-mail: bakirova1974@mail.ru, $~^3$e-mail: ruteshova1@gmail.com%
\begingroup
\renewcommand{\thefootnote}{}
\footnotetext{%
\includegraphics[height=7.0mm]{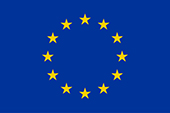}\,
This paper is supported by the European Union's Horizon 2020
research and innovation programme under the Marie
Sklodowska-Curie grant agreement ID: 873071, project SOMPATY
(Spectral Optimization: From Mathematics to Physics and Advanced
Technology).%
}
\endgroup

\vspace{5mm}

{\bf Abstract}

{\footnotesize In the present paper, we study a multipoint boundary value problem for a system of Fredholm integro-differenial equations by the method of parameterization.The case of a degenerate kernel is studied separately, for which we obtain well-posedness conditions and propose some algorithms to find approximate and numerical solutions of the problem. We then establish necessary and sufficient conditions for the well-posedness of the multipoint problem for a system of Fredholm integro-differential equations and develop some algorithms for finding its approximate solutions. These algorithms are based on the solutions of an approximating problem for the system of integro-differential equations with degenerate kernel.}

\vskip0.3cm



\vskip0.5cm

{\bf 1. Introduction}

\vskip0.2cm

Various types of multipoint problems for differential and integro-differential equations have been studied by many researchers, see [1-4, 7-11, 20-25]. A number of methods have been applied to solve these problems, e.g., methods of qualitative theory of differential equations, the method of Green's functions, the method of upper and lower solutions, numerical-analytical methods. However, the problem of establishing effective criteria for the unique solvability of multipoint problems for integro-differential equations, as well as developing algorithms for finding their approximate and numerical solutions, still remains open.

One of constructive methods of investigation and solving boundary value problems for ordinary differential equations and integro-differential equations is the method of parameterization proposed by Dzhumabaev [12]. This method was originally developed for studying and solving boundary value problems for systems of ordinary differential equations. In [12], coefficient criteria were established for the unique solvability of linear boundary value problems. An algorithm for finding their approximate solutions was developed. The method of parameterization was later extended to linear multipoint boundary value problems [20-21], for which necessary and sufficient conditions were obtained for the unique solvability in terms of initial data and an algorithm for finding their approximate solutions was proposed. In [13-15, 19], the method of parameterization was applied to two-point boundary value problems for Fredholm integro-differential equations to establish criteria for their solvability and unique solvability. For these problems, based on the method of parameterization and a new concept of general solution, novel algorithms for approximate and numerical solutions were developed, see [16-18]. The results obtained in above-mentioned papers were used to investigate a multipoint boundary value problem for loaded differential equations [4] and a boundary value problem with a parameter for Fredholm integro-differential equations [3].

Consider the multipoint boundary value problem for the system of integro-differential equations
$$\frac{dx}{dt}=A(t)x + \int ^{T}_{0} K(t,\tau)x(\tau)d\tau  + f(t), \qquad x \in R^n,
\quad t\in (0,T), \eqno(1.1)$$
$$\sum \limits ^m_{i=0}B_i x(t_i) = d, \qquad d\in R^{n}.  \eqno (1.2)$$
Here $x(t)=col(x_1(t), x_2(t),..., x_n(t))$ is an unknown function, $(n\times n)$ matrix  $A(t)$ and $n$-vector  $f(t)$ are continuous on $[0,T]$,  $(n\times n)$ matrix $K(t,\tau)$ is continuous on $[0,T]\times [0,T]$, $B_i$ are constant $(n\times n)$ matrices, $0=t_0 < t_1 < t_2 < ... < t_{m-1} < t_m = T$,
$\|x\|=\max
\limits_{i=\overline{1,n}}|x_i|$.

A solution to multipoint problem (1.1), (1.2) is a function $x^{\ast}(t): [0,T]\to \mathbb{R}^n$ that is continuous on $[0,T]$, continuously differentiable on $(0,T)$ and satisfies integro-differential equations (1.1) and multipoint condition (1.2).

The aim of the present paper is to obtain criteria for the unique solvability of problem (1.1), (1.2) and develop algorithms for finding its approximate solutions. To this end, the parameterization method is used. The interval $[0,T]$ is partitioned and additional parameters are introduced as the values of the desired solution at the left endpoints of the partition subintervals. When applying the method of parameterization to problem (1.1), (1.2), some intermediate problems occur, so called special Cauchy problems for integro-differential equations with parameters. The questions of solvability and unique solvability of such problems were thoroughly investigated in [13-19].

Section 2 is devoted to the study of Fredholm integro-differential equations with degenerate kernel. We divide $[0,T]$ into $m$ parts and introduce additional parameters as the values of the desired solution at the left endpoints $t=t_i$, $i=\overline{0,m-1},$ of the subintervals. The unique solvability of a special Cauchy problem for the $\Delta_m$ partition is equivalent to the invertibility of a matrix $I-G(\Delta_m)$ constructed through a fundamental matrix of the differential part and the matrices of the integral kernel. The $\Delta_m$ partition is called regular if the matrix $I-G(\Delta_m)$ is invertible (see [16]).

For a regular $\Delta_m$ partition, a system of linear algebraic equations in the parameters introduced is constructed using $[I-G(\Delta_m)]^{-1},$ the multipoint condition (1.2), and the continuity conditions at the interior partition points $t=t_i$, $i=\overline{1,m-1}.$ It is shown that the invertibility of the matrix of the system constructed is equivalent to the unique solvability of the multipoint boundary value problem under consideration.

In Section 3, we develop the algorithms for finding a solution to a multipoint boundary value problem for the integro-differential equation with degenerate kernel. For a chosen $\Delta_m$ partition, the matrix $G(\Delta_m)$ is calculated. If there is an inverse of $I-G(\Delta_m)$, then we construct a system of linear algebraic equations. The elements of $G(\Delta_m),$
the coefficients and right-hand side of the system are determined by the solutions of the Cauchy problems for ordinary differential equations and the values of the definite integrals of some functions over the partition subintervals. By solving the system of algebraic equations, we determine the values of the solution at the left endpoints of the subintervals. Next, using the values obtained and the data of the integro-differential equation we compose a function
$\mathcal{F}^{\ast}(t)$ that is continuous on $[0,T]$. Solving the Cauchy problems for ordinary differential equations with the right-hand side $\mathcal{F}^{\ast}(t)$, we get the values of the desired solution at the remaining points of the interval $[0,T]$.

If a fundamental matrix of the differential part is found explicitly and the integrals are evaluated exactly, then the algorithm allows us to find a closed-form solution as well. As is known, it is usually impossible to explicitly find a fundamental matrix for a system of ordinary differential equations with variable coefficients, and, in general, only approximate values of definite integrals can be obtained. For this reason, in this section we propose a numerical implementation of the algorithm. The Cauchy problems for ordinary differential equations on the subintervals are solved by the fourth-order Runge-Kutta method; the integrals are calculated by the Simpson formula. It should be noted that the elements of the matrix $G(\Delta_m),$ the coefficients and the right-hand side of the system of algebraic equations in parameters can be evaluated by parallel computing on the partition subintervals.

In Section 4, a multipoint boundary value problem is considered for a Fredholm integro-differential equation when the integral kernel is not degenerate. We approximate the kernel by the degenerate one and then use the results obtained in Section 2.  At each step of the process, a multipoint boundary value problem for the integro-differential equation with degenerate kernel is solved. We establish sufficient conditions for the convergence of the iterative process to a solution of the multipoint boundary value problem for the Fredholm integro-differential equation with non-degenerate kernel. The accuracy of the approximate solution depends on the choice of the approximating kernel and the number of iteration steps. The necessary and sufficient conditions for the well-posedness of the multipoint problem (1.1), (1.2) are obtained in terms of properties of solutions to approximating problems for integro-differential equations with degenerate kernels.

\vskip0.5cm

{\bf 2. The well-posedness of multipoint problems for Fredholm integro-differential equations with degenerate kernel}

\vskip0.2cm

Consider the integro-differential equation
$$ \frac{dx}{dt} =
A(t)x + \sum\limits_{j=1}^{k} \int ^T_0 \varphi_j(t)
\psi_{j}(\tau)x(\tau)d\tau +f(t), \quad t\in (0,T), \qquad x\in R^n, \eqno
(2.1)$$ subject to the multipoint condition
$$ \sum \limits ^m_{i=0}B_i x(t_i) = d,  \qquad  d \in R^n, \eqno (2.2)$$
where the matrices $A(t),$ $\varphi_j(t),$ $\psi_j(\tau),$
$j=\overline{1,k},$ and the vector $f(t)$ are continuous on $[0,T],$
$\|x\|=\max\limits_{i=\overline{1,n}}|x_i|.$

 The interval $[0,T)$ is divided into $m$ parts by the points $t_0=0 <t_1 <\ldots< t_m=T$, and the partition $\displaystyle{[0,T)=\bigcup\limits_{r=1}^m[t_{r-1},t_r)}$ is denoted by $\Delta_m$.  The case of no partitioning the interval $[0,T]$ is denoted by $\Delta_1$.

We introduce the following spaces: $C([0,T],R^{n})$  is the space of continuous functions $x:[0,T]\rightarrow R^{n}$ with the norm
$\|x\|_1= \max \limits _{t\in[0,T]}\|x(t)\|;$

$C([0,T],\Delta_N,R^{nm})$ is the space of function systems
$x[t]=(x_1(t),x_2(t),\ldots,x_m(t)),$ where functions $x_r:[t_{r-1},t_r)\rightarrow R^{n}$ are continuous and have finite left-handed limits $\lim\limits _{t\rightarrow t_r-0}x_r(t)$
for all $r=\overline{1,m},$ with the norm
$\|x[\cdot]\|_2=\max\limits_{r=\overline{1,m}}\sup\limits_{t\in[t_{r-1},t_r)}
\|x_r(t)\|.$

Let $x_r(t)$ be the restriction of the function $x(t)$ to the $r$th
subinterval $[t_{r-1}, t_r),$ i.e. $x_r(t)=x(t),$ $t\in[t_{r-1},
t_r),$ $r=\overline{1,m}.$

We introduce additional parameters $\lambda _r=x_r(t_{r-1})$ and make the substitution $ x_r(t) = u_r(t) +\lambda_r$ on each $r$th subinterval. The multipoint problem (2.1), (2.2) is then reduced to the following problem with parameters:
$$ \frac{du_r}{dt} =
A(t)(u_r+\lambda
_r)+\sum\limits_{s=1}^{m}\sum\limits_{j=1}^{k}\int ^{t_s}_{t_{s-1}}\varphi_j(t)\psi_j(\tau)[u_r(\tau) + \lambda_r] d\tau +
f(t), \,\, t \in (t_{r-1},t_r),\,\, r= \overline{1,m},\eqno (2.3)
$$
$$ u_r(t_{r-1})= 0, \qquad r= \overline{1,m}, \eqno (2.4) $$
$$ \sum \limits ^{m-1}_{i=0}B_i \lambda_{i+1} + B_m\lambda _m +  B_m \lim \limits _{t\rightarrow T-0} u_m(t) = d, \eqno (2.5) $$
$$ \lambda _{p} + \lim \limits_{t\to t_p-0} u_{p}(t) - \lambda_{p+1}=0, \qquad
p=\overline{1,m-1}, \eqno (2.6) $$ where (2.6) are the continuity conditions for the solution at the interior points of the partition $\Delta_m.$ Note that conditions (2.6) and integro-differential equations (2.3) ensure the continuity of the derivative of the solution at those points.

If $x^{\ast}(t)$ is a solution to multipoint problem (2.1),(2.2), then the pair
 $(\lambda^{\ast},u^{\ast}[t])$ with elements
$\lambda^{\ast}=(\lambda_1^{\ast},\lambda_2^{\ast},\ldots,\lambda_m^{\ast})\in
R^{nm},$ $u^{\ast}[t] =\bigl{(}u_1^{\ast}(t),u_2^{\ast}(t),\ldots,
u_m^{\ast}(t)\bigr{)}\in C([0,T],\Delta_m,R^{nm}),$ where
$\lambda_r^{\ast}=x^{\ast}(t_{r-1}),$
$u_r^{\ast}(t)=x^{\ast}(t)-x^{\ast}(t_{r-1}),$  $[t_{r-1},t_r),$
$r=\overline{1,m}$, is a solution to the problem with parameters (2.3)-(2.6).
Vice versa, if a pair $(\widetilde{\lambda}, \widetilde{u}[t])$
with elements
$\widetilde{\lambda}=(\widetilde{\lambda}_1,\widetilde{\lambda}_2,\ldots,\widetilde
{\lambda}_m)\in R^{nm},$ $\widetilde{u}[t]
=\bigl{(}\widetilde{u}_1(t),\widetilde{u}_2(t),\ldots,
\widetilde{u}_m(t)\bigr{)}\in C([0,T],\Delta_m,R^{nm}),$ is a solution to problem with parameters (2.3)-(2.6), then the function
$\widetilde{x}(t)$ defined as
$\widetilde{x}(t)=\widetilde{\lambda}_r+\widetilde{u}_r(t),$
$t\in[t_{r-1},t_r),$ $r=\overline{1,m},$
$\widetilde{x}(T)=\widetilde{\lambda}_m+\lim\limits_{t\rightarrow
T-0}\widetilde{u}_m(t),$ is a solution to the original problem (2.1),(2.2).

If $X_r(t)$ is a fundamental matrix of the differential equation
$\displaystyle{ \frac{dx_r}{dt} = A(t)x_r}$ on $[t_{r-1},t_r]$, then the special Cauchy problem for the system of integro-differential equations with parameters (2.3), (2.4) is reduced to the equivalent system of integral equations
$$ u_r(t) =
X_r(t)\int _{t_{r-1}}^{t}X_r^{-1}(\tau)A(\tau)d \tau
\lambda_r+ X_r(t)\int _{t_{r-1}}^{t}X_r^{-1}(\tau)\sum\limits_{s=1}^{m}
\sum\limits_{j=1}^{k} \int _{t_{s-1}}^{t_s}
\varphi_{j}(\tau)\psi_j(\tau_1)u_s(\tau_1)d\tau_1 d\tau+$$
$$+ X_r(t)\int _{t_{r-1}}^{t}X_r^{-1}(\tau)\sum\limits_{s=1}^{m}
\sum\limits_{j=1}^{k} \int _{t_{s-1}}^{t_s}
\varphi_{j}(\tau)\psi_j(\tau_1)d\tau_1 d\tau \lambda_s +
X_r(t)\int_{t_{r-1}}^{t}X_r^{-1}(\tau)f(\tau) d\tau, \eqno (2.7) $$

$ \qquad  \qquad
t \in [t_{r-1},t_r), \qquad r=\overline{1,m}.$

 Setting $\displaystyle \mu_j =\sum\limits_{s=1}^{m}
\int _{t_{s-1}}^{t_s}\psi_j(\tau)u_s(\tau)d\tau$, rewrite (2.7) in the following way:
$$ u_r(t)=\sum\limits_{j=1}^{k}X_r(t)\int _{t_{r-1}}^{t}X_r^{-1}(\tau)
\varphi_{j}(\tau)d\tau \mu_j +
X_r(t)\int _{t_{r-1}}^{t}X_r^{-1}(\tau)\biggl[
A(\tau) \lambda_r+ $$ $$ + \sum\limits_{j=1}^{k}
\varphi_{j}(\tau)\sum\limits_{s=1}^{m} \int _{t_{s-1}}^{t_s}
\psi_j(\tau_1)d\tau_1 \lambda_s + f(\tau)\biggr]d\tau, \qquad t \in [t_{r-1},t_r),
\qquad r=\overline{1,m}. \eqno (2.8) $$

Multiplying both sides of (2.8) by $\psi_p(t)$, integrating them over $[t_{r-1},t_r]$,
and summing up with respect to $r$, we get the following system of linear algebraic equations in  $\mu=(\mu_1,\ldots,\mu_k)\in R^{nk}$:
$$\mu_p=\sum\limits_{l=1}^k G_{p,l}(\Delta_m)\mu_l + \sum\limits^m_{r=1}V_{p,r}(\Delta_m)\lambda_r+
g_p(f,\Delta_m), \quad   p=\overline{1,k}, \eqno (2.9)$$ with $(n\times n)$ matrices
$$G_{p,l}(\Delta_m)=
\sum\limits_{r=1}^{m}\int_{t_{r-1}}^{t_r} \psi_p(\tau)
X_r(\tau)\int_{t_{r-1}}^{\tau} X_r^{-1}(\tau_1)\varphi_l(\tau_1) d\tau_1
d\tau, \eqno (2.10) $$
$$V_{p,r}(\Delta_m)=
\int _{t_{r-1}}^{t_r} \psi_p(\tau)
X_r(\tau)\int _{t_{r-1}}^{\tau} X^{-1}_r(\tau_1)A(\tau_1) d\tau_1
d\tau+$$ $$+\sum \limits_{s=1}^{m}\sum\limits_{j=1}^{k}\int_{t_{s-1}}^{t_s}
\psi_p(\tau) X_s(\tau)\int_{t_{s-1}}^{\tau}
X_s^{-1}(\tau_1)\varphi_j(\tau_1) d\tau_1 d\tau   \int _{t_{r-1}}^{t_r}\psi_j(\tau_2)d\tau_2,  \eqno (2.11)
$$ and vectors of dimension $n$
$$g_{p}(f,\Delta_m)=
\sum\limits_{r=1}^{m}\int_{t_{r-1}}^{t_r} \psi_p(\tau)
X_r(\tau)\int_{t_{r-1}}^{\tau} X_r^{-1}(\tau_1)f(\tau_1)d\tau_1 d\tau,
\quad p=\overline{1,k}, \quad j=\overline{1,k}. \eqno (2.12)
$$
Using the matrices $G_{p,l}(\Delta_m)$ and $V_{p,r}(\Delta_m)$, we construct the matrices  $G(\Delta_m)=(G_{p,l}(\Delta_m)),$
$p,l=\overline{1,k},$ and $V(\Delta_m)=(V_{p,r}(\Delta_m)),$
$p=\overline{1,k},$ $r=\overline{1,m}$. Then, system (2.9) can be rewritten in the form
$$[I-G(\Delta_m)]\mu=V(\Delta_m)\lambda+g(f,\Delta_m), \eqno (2.13)$$
where $I$ is the identity matrix of dimension $nk,$
$g(f,\Delta_m)=(g_1(f,\Delta_m),\ldots,g_k(f,\Delta_m))\in
R^{nk}.$

{\bf Definition 2.1} The partition $\Delta_m$ is said to be regular if the matrix $I-G(\Delta_m)$ is invertible.

Any fundamental matrix of the differential equation
$\displaystyle{ \frac{dx_r}{dt} = A(t)x_r}$ on $[t_{r-1},t_r]$ can be represented as $X_r(t)=X_r^0(t)\cdot C_r,$ where  $X_r^0(t)$ is the normalized fundamental matrix $(X_r^0(t_{r-1})=I)$ and $C_r$ is an arbitrary invertible matrix. Thus, the following equalities hold true
$$G_{p,l}(\Delta_m)=
\sum \limits_{r=1}^{m}\int _{t_{r-1}}^{t_r} \psi_p(\tau)
X_r^0(\tau)C_r\int _{t_{r-1}}^{\tau}[X_r^{0}(\tau_1)C_r]^{-1}\varphi_l(\tau_1)
d\tau_1 d\tau=$$$$= \sum\limits_{r=1}^{m}\int _{t_{r-1}}^{t_r}
\psi_p(\tau) X_r^0(\tau)\int _{t_{r-1}}^{\tau}
[X_r^{0}(\tau_1)]^{-1}\varphi_l(\tau_1) d\tau_1 d\tau,$$ and the regularity of the $\Delta_m$ partition does not depend on the choice of a fundamental matrix for the differential part of the equation.

Let us denote by $\sigma(k,[0,T])$ the set of regular partitions $\Delta_m$ of the interval $[0,T]$ for Eq. (2.1).

{\bf Definition 2.2} The special Cauchy problem (2.3),(2.4) is called uniquely solvable if it has a unique solution for any $\lambda\in R^{nm}$ and $f(t)\in C([0,T],R^n)$ .

The special Cauchy problem (2.3), (2.4) is equivalent to the system of integral equations (2.7). Since the kernel of (2.7) is degenerate, this system, in turn, is equivalent to the system of algebraic equations (2.9) in $\mu=(\mu_1,\ldots,\mu_k)\in R^{nk}.$
Therefore, the special Cauchy problem is uniquely solvable if and only if the $\Delta_m$ partition, generating this problem, is regular.

Since the special Cauchy problem is uniquely solvable for a partition with a sufficiently small step size $h>0$ (see [13, p.1152]), the set $\sigma(k,[0,T])$ is not empty.

Let us take a partition $\Delta_m\in \sigma(k,[0,T])$ and represent the matrix  $[I-G(\Delta_m)]^{-1}$ in the form $[I-G(\Delta_m)]^{-1}=\Big(M_{j,p}(\Delta_m)\Big),$
$j, p=\overline{1,k},$ where $M_{j,p}(\Delta_m)$ are square matrices of dimension $n.$

Then, taking into account (2.13), we can determine the elements of the vector $\mu\in R^{nk}$
from the equalities
$$\mu_j=\sum\limits_{i=1}^{m}\Big (\sum\limits_{p=1}^{k} M_{j,p}(\Delta_m)V_{p,i}(\Delta_m)\Big)\lambda_i +
\sum\limits^k_{p=1}M_{j,p}(\Delta_m)g_p(f,\Delta_m), \quad
j=\overline{1,k}. \eqno (2.14)$$ In (2.8), by replacing $\mu_j$ with the right-hand side of (2.14) we get the representation of the functions $u_r(t)$ through $\lambda_i,$
$i=\overline{1,m}:$
$$ u_r(t)=\sum\limits_{i=1}^{m}\biggl\{ \sum\limits_{j=1}^{k}X_r(t)\int _{t_{r-1}}^{t}X_r^{-1}(\tau)
\varphi_{j}(\tau)d\tau \biggl [\sum\limits_{p=1}^{k} M_{j,p}(\Delta_m)V_{p,i}(\Delta_m) + \int _{t_{i-1}}^{t_i}
\psi_j(\tau_1)d\tau_1 \biggr ] \biggr\}\lambda_i
 + $$
 $$ + X_r(t)\int _{t_{r-1}}^{t}X_r^{-1}(\tau)
A(\tau)d\tau \lambda_r+ $$ $$ +  X_r(t)\int _{t_{r-1}}^{t}X_r^{-1}(\tau) \biggl [
\sum\limits_{j=1}^{k} \varphi_{j}(\tau)\sum\limits^k_{p=1}M_{j,p}(\Delta_m)g_p(f,\Delta_m) +  f(\tau)\biggr]d\tau, \quad t \in [t_{r-1},t_r),
\ r=\overline{1,m}.
 \eqno (2.15) $$
We will use the following notation:
$$D_{r,i}(\Delta_m)=\sum\limits_{j=1}^{k}X_r(t_r)\int _{t_{r-1}}^{t}X_r^{-1}(\tau)
\varphi_{j}(\tau)d\tau \biggl [\sum\limits_{p=1}^{k} M_{j,p}(\Delta_m)V_{p,i}(\Delta_m) + $$ $$ + \int _{t_{i-1}}^{t_i}
\psi_j(\tau_1)d\tau_1 \biggr ], \qquad i\neq r,  \qquad r,j=\overline{1,m},\eqno (2.16)$$
$$D_{r,r}(\Delta_m)= \sum\limits_{j=1}^{k}X_r(t_r)\int _{t_{r-1}}^{t}X_r^{-1}(\tau)
\varphi_{j}(\tau)d\tau \biggl [\sum\limits_{p=1}^{k} M_{j,p}(\Delta_m)V_{p,r}(\Delta_m) + \int _{t_{r-1}}^{t_r}
\psi_j(\tau_1)d\tau_1 \biggr ] + $$ $$ + X_r(t_r)\int _{t_{r-1}}^{t}X_r^{-1}(\tau)
A(\tau)d\tau, \eqno (2.17)$$
$$F_{r}(\Delta_m)=X_r(t)\int _{t_{r-1}}^{t}X_r^{-1}(\tau) \biggl [
\sum\limits_{j=1}^{k} \varphi_{j}(\tau)\sum\limits^k_{p=1}M_{j,p}(\Delta_m)g_p(f,\Delta_m) +  f(\tau)\biggr]d\tau, \quad r=\overline{1,m}.\eqno (2.18)$$

From (2.15), we find the limits
$$\lim\limits_{t\rightarrow t_r-0}u_r(t)=\sum\limits_{i=1}^m D_{r,i}(\Delta_m)\lambda_i + F_r(\Delta_m).
\eqno (2.19)$$

Substituting the right-hand side of (2.19) into condition (2.5) and continuity conditions (2.6),
we get the following system of linear algebraic equations in parameters $\lambda_r,$ $r=\overline{1,m}$:
$$\sum \limits ^{m-2}_{i=0}[B_i+ B_m D_{r,i+1}(\Delta_m)]  \lambda_{i+1} + [B_{m-1} + B_m + B_m D_{r,m}(\Delta_m)]\lambda _m  = d - B_m F_m(\Delta_m), \eqno (2.20) $$
$$[I+D_{p,p}(\Delta_m)]\lambda_p-[I-D_{p,p+1}(\Delta_m)]\lambda_{p+1}+\sum\limits_{{}^{\,\,\,\,\,\,\,\,\,\,i=1}_{i \neq p, \, i \neq p+1}}^{m}D_{p,i}(\Delta_m)\lambda_i=-F_p(\Delta_m),
\quad p=\overline{1,m-1}. \eqno(2.21)$$

Let $Q_{\ast}(\Delta_m)$ be the matrix corresponding to the left-hand side of system (2.20), (2.21). Then the system (2.20), (2.21) can be written in the form
 $$Q_{\ast}(\Delta_m)\lambda=-F_{\ast}(\Delta_m), \quad \lambda\in R^{nm},  \eqno (2.22)$$
where
$F_{\ast}(\Delta_m)=\Big(-d+B_m F_{m}(\Delta_m),F_{1}(\Delta_m),\ldots,F_{m-1}(\Delta_m)\Big)\in
R^{nm}.$

{\bf Lemma 2.1}  {\it The following statements hold true for $\Delta_m \in \sigma(k,[0,T])$:\\
(a) the vector
$\lambda^{\ast}=(\lambda_1^{\ast},\lambda_2^{\ast},\ldots,\lambda_m^{\ast})\in
R^{nm},$ composed of the values of a solution $x^{\ast}(t)$ to problem (2.1),(2.2) at the partition points, $\lambda_r^{\ast}=x^{\ast}(t_{r-1}),$ $r=\overline{1,m}$, satisfies system (2.22); \\
(b) if
$\widetilde{\lambda}=(\widetilde{\lambda}_1,\widetilde{\lambda}_2,\ldots,\widetilde{\lambda}_m)\in
R^{nm}$ is a solution to system (2.22) and a function system  $\widetilde{u}[t]=(\widetilde{u}_{1}(t),\widetilde{u}_{2}(t),...,\widetilde{u}_{m}(t))$
is a solution to the special Cauchy problem  (2.3), (2.4) with
$\lambda_r=\widetilde{\lambda}_r,$ $r=\overline{1,m},$ then the function  $\widetilde{x}(t),$  defined as
$\widetilde{x}(t)=\widetilde{\lambda}_r + \widetilde{u}_r(t),$ $t\in
[t_{r-1},t_r),$  $r=\overline{1,m},$
$\widetilde{x}(T)=\widetilde{\lambda}_m + \lim\limits_{t\rightarrow
T-0}\widetilde{u}_m(t),$ is a solution to problem (2.1),(2.2).}

The proof of Lemma 2.1 repeats the proof of Lemma 1 in [6, p. 1155] with minor changes.

We will use the following notation:
$$\alpha=\max\limits_{t\in[0,T]}\|A(t)\|,  \quad
\overline{\omega}=\max\limits_{r=\overline{1,m}}(t_r-t_{r-1}),$$
$$\overline{\varphi}(k)=\max\limits_{r=\overline{1,m}}\int\limits_{t_{r-1}}^{t_r}
\sum\limits_{j=1}^{k} \|\varphi_j(t)\|dt, \quad
\overline{\psi}(T)=
\max\limits_{p=\overline{1,k}}\int\limits_{0}^{T}
\|\psi_p(t)\|dt.$$

{\bf Theorem 2.1} {\it Let $\Delta_m \in \sigma(k,[0,T])$ and the matrix  ${Q}_{\ast}(\Delta_m): R^{nm}\to R^{nm}$ be invertible. Then problem (2.1),(2.2) has a unique solution $x^{\ast}(t)$ for any
$f(t)\in C([0,T],R^{n}),$ $d \in R^{n},$ and the following estimate holds:
$$\|x^{\ast}\|_1\leq\mathcal{N}(k,\Delta_m)\max (\|d\|,\|f\|_1), \eqno (2.23)$$
where
$$\mathcal{N}(k,\Delta_m)=e^{\alpha\overline{\omega}}\Big\{\overline{\varphi}(k)
\Big[\|[I-G(\Delta_m)]^{-1}\| \cdot\overline{\psi}(T)
\Big(e^{\alpha\overline{\omega}}-1+e^{\alpha\overline{\omega}}
\cdot\overline{\varphi}(k)\cdot\overline{\psi}(T)\Big)+\overline{\psi}(T)\Big]+1\Big\}
\times$$
$$\times\gamma_{\ast}(\Delta_m)
(1+\|C\|)\max\Big\{1,\overline{\omega}e^{\alpha\overline{\omega}}\Big[1+
e^{\alpha\overline{\omega}}\cdot\overline{\varphi}(k)\cdot\|[I-G(\Delta_m)]^{-1}\|\cdot
\overline{\psi}(T)\Big\}+$$$$+
e^{\alpha\overline{\omega}}\overline{\omega}\Big[\overline{\varphi}(k)\cdot\|[I-G(\Delta_m)]^{-1}\|\cdot
\overline{\psi}(T)\cdot e^{\alpha\overline{\omega}}+1\Big]. \eqno
(2.24)$$}
 {\bf Proof.} Take a partition $\Delta_m \in  \sigma(k,[0,T])$. Let $f(t)\in C([0,T],R^{n})$  and $d\in R^{n}$.
 Since the matrix ${Q}_{\ast}(\Delta_m)$ is invertible, we can find the unique solution to the system of linear algebraic equations (2.22):
$$\lambda^{\ast}=-[Q_{\ast}(\Delta_m)]^{-1}F_{\ast}(\Delta_m).  $$
By solving the special Cauchy problem (2.3), (2.4) with $\lambda=\lambda^{\ast}$, we get the function system
$u^{\ast}[t]=(u^{\ast}_1(t),u^{\ast}_2(t),\ldots,u^{\ast}_m(t)).$
It follows from the regularity of the $\Delta_m$ partition that there exists a unique function system $u^{\ast}[t]$ with the elements
$u^{\ast}_r(t)$ that are determined from the right-hand side of representation (2.15) with
$\lambda=\lambda^{\ast}=(\lambda^{\ast}_1,\lambda^{\ast}_2,\ldots,\lambda^{\ast}_m)\in
R^{nm}.$ Then, by Lemma 2.1, the function  $x^{\ast}(t)$
defined as
$x^{\ast}(t)=\lambda_r^{\ast} + u_r^{\ast}(t),$ $\ t\in[t_{r-1},t_r),$
$ r=\overline{1,m},$
$ x^{\ast}(T)=\lambda_m^{\ast} + \lim\limits_{t\rightarrow
T-0}u_m^{\ast}(t),$ is a solution to problem (2.1), (2.2).
The uniqueness of a solution can be proved by contradiction.

Let us verify the validity of the estimate (2.23).

Using the equalities
$$ X_r(t)\int _{t_{r-1}}^{t}X_r^{-1}(\tau)P(\tau)d\tau= \int _{t_{r-1}}^{t}P(\tau_1)d\tau_1+
\int _{t_{r-1}}^{t}A(\tau_1)\int _{t_{r-1}}^{\tau_1}P(\tau_2)d\tau_2d\tau_1+$$ $$+
\int _{t_{r-1}}^{t}A(\tau_1)\int _{t_{r-1}}^{\tau_1}A(\tau_2)
\int _{t_{r-1}}^{\tau_2}P(\tau_3)d\tau_3d\tau_2d\tau_1+\ldots,
\quad t\in [t_{r-1},t_r],$$ we get the estimates
$$\|X_r(t_r)\int _{t_{r-1}}^{t_r}X_r^{-1}(\tau)\varphi_j(\tau)d\tau\|\leq e^{\alpha(t_r-t_{r-1})}\int _{t_{r-1}}^{t_r}
\|\varphi_j(t)\|dt, \quad   r=\overline{1,m}. \eqno (2.25)$$

It follows from (2.12), (2.18), (2.25) that
$$\|g_{p}(f,\Delta_m)\|\leq
\sum\limits_{r=1}^{m}\int _{t_{r-1}}^{t_r}\|
\psi_p(\tau)\|\cdot \|X_r(\tau)\int _{t_{r-1}}^{\tau}
X_r^{-1}(\tau_1)f(\tau_1)d\tau_1\|d\tau\leq $$ $$ \leq
\sum\limits_{r=1}^{m}\int _{t_{r-1}}^{t_r}\|
\psi_p(\tau)\|d\tau\cdot
e^{\alpha\overline{\omega}}\cdot\overline{\omega}\cdot\|f\|_1=\int _{0}^{T}\|
\psi_p(t)\|dt\cdot
e^{\alpha\overline{\omega}}\cdot\overline{\omega}\cdot\|f\|_1,
\quad p=\overline{1,k}, \eqno (2.26)$$
$$\|F_{r}(\Delta_m)\|\leq e^{\alpha(t_r-t_{r-1})}
\sum\limits_{j=1}^{k}\int _{t_{r-1}}^{t_r}\| \psi_j(t)\|dt
\|[I-G(\Delta_m)]^{-1}\|\max\limits_{p=\overline{1,k}}\|g_p(f,\Delta_m)\|+
e^{\alpha(t_r-t_{r-1})} \omega_r \|f\|_1.$$

By using
$$\|F_{\ast}(\Delta_m)\|\leq (1+\|B_m\|) \max
\Bigl(\|d\|,\max\limits_{r=\overline{1,m}}\|F_r(\Delta_m)\|\Bigr),$$
 and taking into account (2.16),(2.25), and (2.26), we get
$$\|F_{\ast}(\Delta_m)\|\leq
(1+\|B_m\|)\max\Big\{1, \overline{\omega}
e^{\alpha\overline{\omega}}\Big[1+  $$ $$+ e^{\alpha\overline{\omega}}
\max\limits_{r=\overline{1,m}}\int _{t_{r-1}}^{t_r}\sum\limits_{j=1}^{k}
\|\varphi_j(t)\|dt\|[I-G(\Delta_m)]^{-1}\|
\max\limits_{p=\overline{1,k}}\int _{0}^{T}\|
\psi_p(t)\|dt\Big]\Big\}\max(\|d\|,\|f\|_1). \eqno (2.27)$$

The inequalities (2.22), (2.27) and the invertibility of
$Q_{\ast}(\Delta_m)$  yield the following estimate:
$$\|\lambda^{\ast}\|
\leq\|[Q_{\ast}(\Delta_m)]^{-1}\|\|F_{\ast}
(\Delta_m)\|\leq\gamma_{\ast}(\Delta_m)(1+\|B_m\|)\max\Big\{1,
\overline{\omega}
e^{\alpha\overline{\omega}}\Big[1+ $$ $$ + e^{\alpha\overline{\omega}}
\max\limits_{r=\overline{1,m}}\int _{t_{r-1}}^{t_r}\sum\limits_{j=1}^{k}
\|\varphi_j(t)\|dt \|[I-G(\Delta_m)]^{-1}\|
\max\limits_{p=\overline{1,k}}\int _{0}^{T}\|
\psi_p(t)\|dt\Big\}\max(\|d\|,\|f\|_1). \eqno (2.28)$$
By (2.15) and (2.11), we get
$$\|u^{\ast}[\cdot]\|_2\leq \biggl\{e^{\alpha\overline{\omega}}
\max\limits_{r=\overline{1,m}}\int _{t_{r-1}}^{t_r}\sum\limits_{j=1}^{k}\|\varphi_j(t)\|dt\|[I-G(\Delta_m)]^{-1}\|
\max\limits_{p=\overline{1,k}}\int _{0}^{T}\|
\psi_p(t)\|dt \Big(e^{\alpha\overline{\omega}}-1+ $$ $$ + e^{\alpha\overline{\omega}}
\max\limits_{r=\overline{1,m}}\int _{t_{r-1}}^{t_r}
\max\limits_{j=\overline{1,k}} \|\varphi_j(t)\|dt
\max\limits_{p=\overline{1,k}}\int _{0}^{T}\|
\psi_p(t)\|dt\Big) +
\max\limits_{p=\overline{1,k}}\int _{0}^{T}\|
\psi_p(t)\|dt\Big]+(e^{\alpha\overline{\omega}}-1)\biggr\}\|\lambda\|+$$ $$+
e^{\alpha\overline{\omega}}\overline{\omega}
\Big[\max\limits_{r=\overline{1,m}}\int _{t_{r-1}}^{t_r}\sum\limits_{j=1}^{k}
\|\varphi_j(t)\|dt\|[I-G(\Delta_m)]^{-1}\|
\max\limits_{p=\overline{1,k}}\int _{0}^{T}\|
\psi_p(t)\|dt\cdot e^{\alpha \overline{\omega}}+1\Big]\cdot\|f\|_1.
\eqno (2.29)$$

Finally, by using (2.28), (2.29) and
$\|x^{\ast}\|_1\leq\|\lambda^{\ast}\|+\|u^{\ast}[\cdot]\|_2$, we arrive at the estimate (2.23).

Theorem 2.1 is proved.

{\bf Definition 2.3} Problem (2.1), (2.2) is said to be well-posed if, for any pair $(f(t),d),$ with $f(t)\in C([0,T],R^n)$ and $d\in R^n$, it has a unique solution $x(t)$, and the estimate
$$\|x\|_1\leq K\max(\|f\|_1,\|d\|)$$ holds, where  $K$ is a constant independent of $f(t)$ of $d.$

{\bf Theorem 2.2} { \it Problem (2.1),(2.2) is well-posed if and only if the matrix
$Q_{\ast}(\Delta_m):$ $R^{nm}\rightarrow R^{nm}$ is invertible for any partition $\Delta_m \in \sigma(k,[0,T])$.}

{\bf  Proof.} For a fixed $k$ and $\Delta_m \in \sigma(k,[0,T])$ the number  $\mathcal{N}(k,\Delta_m)$, defined by (2.24), does not depend on  $f(t)$ and $d.$ Thus the sufficiency of the conditions of Theorem 2.2 for the well-posedness of problem (2.1), (2.2) follows from Theorem 2.1.

{\it Necessity}. Let problem (2.1), (2.2) be well-posed and $\Delta_m \in \sigma(k,[0,T]).$  Suppose to the contrary that the matrix $Q_{\ast}(\Delta_m):$ $R^{nm}\rightarrow R^{nm}$ is not invertible. This is possible only if the homogeneous system of equtions
$$Q_{\ast}(\Delta_N)\lambda=0,  \quad \lambda\in R^{nN},  \eqno (2.30)$$
has a nonzero solution.
Assuming that
$\widetilde{\lambda}=(\widetilde{\lambda}_1,\widetilde{\lambda}_2,\ldots,\widetilde{\lambda}_m)$
is a nonzero solution (i.e. $\|\widetilde{\lambda}\|\neq 0$) to system (2.30), consider the homogeneous problem (2.1), (2.2) with $f(t)=0$ and  $d=0.$  For this problem, system (2.22) coincides with (2.30). Then, by Lemma 2.1, the function $\widetilde{x}(t)$ defined as
$\widetilde{x}(t)=\widetilde{\lambda}_r + \widetilde{u}_r(t),$ $\ t\in
[t_{r-1},t_r),$ $r=\overline{1,m},$
$ \widetilde{x}(T)=\widetilde{\lambda}_m +\lim\limits_{t\rightarrow
T-0}\widetilde{u}_m(t),$ is a nonzero solution to the homogeneous problem. Here the function system
$\widetilde{u}[t]= (\widetilde{u}(t),\widetilde{u}_2(t),\ldots,
\widetilde{u}_m(t))$ is a solution to the special Cauchy problem
(2.3), (2.4) with $\lambda=\widetilde{\lambda}$ and  $f(t)=0.$
This contradicts the well-posedness of problem (2.1), (2.2).

Theorem 2.2 is proved.

\vskip0.5cm

{\bf 3.} {\bf An algorithm for solving multipoint problems for Fredholm integro-differential equations with degenerate kernel and its numerical implementation}

The Cauchy problems for ordinary differential equations on the subintervals
$$ \frac{dx}{dt} =
A(t)x + P(t), \qquad x(t_{r-1})=0, \quad t\in[t_{r-1},t_r],\quad
r=\overline{1,m}, \eqno (3.1)$$ are an essential part of the algorithm proposed. Here $P(t)$ is an   $(n\times n)$ matrix or $n$ vector that is continuous on $[t_{r-1},t_r],$ $r=\overline{1,m}$.
Hence, a solution to problem (3.1) is a matrix or a vector of dimension $n$.

Let $E_{\ast,r}(A(\cdot),P(\cdot),t)$ denote a solution to the Cauchy problem (3.1). Clearly,
$$E_{\ast,r}(A(\cdot),P(\cdot),t)=X_r(t)\int _{t_{r-1}}^{t}
X^{-1}(\tau)P(\tau)d\tau, \quad  t\in[t_{r-1},t_r], \eqno (3.2)$$
where $X_r(t)$ is a fundamental matrix of differential equation (3.1) on the $r$th subinterval.

The choice of a regular partition is another important part of the algorithm. We can start with $\Delta_1,$ when the interval $[0,T]$ is not partitioned.

I. We divide $[0,T]$ into $m$ parts by the points $t_0=0<t_1<\ldots<t_{m-1}<t_m=T,$ involved in the multipoint condition. The resulting partition we denote by $\Delta_m, m=1,2,\ldots$.

II. By solving $mk$ Cauchy problems for ordinary differential matrix equations $$ \frac{dx}{dt} = A(t)x + \varphi_k(t),
\qquad x(t_{r-1})=0, \quad t\in[t_{r-1},t_r], \eqno (3.3)$$ we obtain the matrix functions
$$E_{\ast,r}(A(\cdot),\varphi_j(\cdot),t), \quad  t\in[t_{r-1},t_r], \quad r=\overline{1,m}, \quad j=\overline{1,k}.  \eqno (3.4)$$

III. We multiply each $(n\times n)$ matrix (3.4) by $(n\times
n)$ matrix $\psi_p(t),$ $p=\overline{1,k},$ and integrate the products over
$[t_{r-1},t_r]:$
$$\widehat{\psi}_{p,r}(\varphi_j)=\int_{t_{r-1}}^{t_r}\psi_p(t)
E_{\ast,r}(A(\cdot),\varphi_j(\cdot),t)dt.  \eqno (3.5)$$  Summing up (3.5) with respect to $r$
and taking into account (2.10), (3.2), we get the $(n\times n)$ matrices
$$G_{p,j}(\Delta_m)=\sum\limits_{r=1}^m\widehat{\psi}_{p,r}(\varphi_j),  \quad p, j =\overline{1,k}. $$
We then construct the  $(nk\times nk)$ matrix
$G(\Delta_m)=(G_{p,j}(\Delta_m)),$  $p, j=\overline{1,k},$ and determine whether the matrix  $[I-G(\Delta_m)]:
R^{nk}\rightarrow R^{nk}$ is invertible. If so, we find its inverse and represent it as  $[I-G(\Delta_m)]^{-
1}=(M_{p,j}(\Delta_m)),$ where $M_{p,j}(\Delta_m))$ are $(n\times
n)$ matrices, $p, j=\overline{1,k}$. We then move on to the next step of the algorithm.

If $[I-G(\Delta_m)]$ is not invertible, i.e. the $\Delta_m$ partition is not regular, then we  take a new partition of $[0,T]$ and start the algorithm again. A simple way to choose a new partition is to take $\Delta_{2m},$ dividing each subinterval $\Delta_m$ in half. We add to the points $t=t_i$ of the multipoint condition the points $ (t_i-t_{i-1})/2$, $i=\overline{1,m}$. Then, redesignating all points as $theta_0 = t_0=0,$ $\theta_1=(t_1-t_{0})/2$, $\theta_2= t_1$, $\theta_3= (t_2-t_{1})/2$, $\theta_4= t_2$, ..., $\theta_{2m-1} = (t_m - t_{m-1})/2$,  $\theta_{2m}=t_m=T$, we again get problem (2.1), (2.2) with multipoint conditions at the points $t=\theta_i$, $i=\overline{0,2m}$.

IV. Solving again the Cauchy problems for ordinary differential equations
$$ \frac{dx}{dt} =
A(t)x + A(t), \quad x(t_{r-1})=0, \quad t\in[t_{r-1},t_r], $$
$$ \frac{dx}{dt} =
A(t)x +f(t), \quad x(t_{r-1})=0, \quad t\in[t_{r-1},t_r],\quad
r=\overline{1,m},$$ we obtain $E_{\ast,r}(A(\cdot),A(\cdot),t)$ and
$E_{\ast,r}(A(\cdot),f(\cdot),t),$  $r=\overline{1,m}.$

V. We evaluate the integrals
$$\widehat{\psi}_{p,r}=\int _{t_{r-1}}^{t_r}\psi_p(t)dt,   \quad
\widehat{\psi}_{p,r}(A)=\int _{t_{r-1}}^{t_r}\psi_p(t)E_{\ast,r}(A(\cdot),A(\cdot),t)dt,$$
$$\widehat{\psi}_{p,r}(f)=\int_{t_{r-1}}^{t_r}\psi_p(t)E_{\ast,r}(A(\cdot),f(\cdot),t)dt.$$
From  (2.11), (2.12), and (3.2) we determine the $(n\times n)$ matrices
$$V_{p,r}(\Delta_m)=\widehat{\psi}_{p,r}(A)+\sum\limits_{i=1}^m\sum\limits_{k=1}^m
\widehat{\psi}_{p,i}(\varphi_j)\cdot\widehat{\psi}_{j,r}$$ and the $n$
vectors
$$g_{p}(f,\Delta_m)=\sum\limits_{r=1}^m\widehat{\psi}_{p,r}(\Delta_m),  \quad
p=\overline{1,k},\quad r=\overline{1,m}.$$

VI. We construct the system of linear algebraic equations in parameters
$$Q_{\ast}(\Delta_m)\lambda = -F_{\ast}(\Delta_m),  \qquad \lambda\in R^{nm}.  \eqno (3.6)$$
The elements of the matrix $Q_{\ast}(\Delta_m)$ and the vector  $F_{\ast}(\Delta_m)=(-d+ B_m F_m(\Delta_m),F_1(\Delta_m), \ldots, $ $ F_{m-1}(\Delta_m))\in
R^{nm}$ are determined by the equalities (2.16), (2.17), (2.18),
where, by (3.2), we replace
$X_r(t_r)\int_{t_{r-1}}^{t_r}
X_r^{-1}(\tau)\varphi_j(\tau)d\tau$ and $X_r(t_r)\int_{t_{r-1}}^{t_r} X_r^{-1}(\tau)f(\tau)d\tau$
with $E_{\ast,r}(A(\cdot),\varphi_j(\cdot),t_r)\ $ and  $ \  E_{\ast,r}(A(\cdot),f(\cdot),t_r),$ respectively.

It follows from Theorem 2.2 that the invertibility of the matrix
$Q_{\ast}(\Delta_m)$ is equivalent to the well-posedness of problem (2.1),(2.2).

By solving system (3.6), we get $\lambda^{\ast}=(\lambda_1^{\ast},\lambda_2^{\ast},\ldots,\lambda_m^{\ast})\in
R^{nm}.$

VII. From the equalities
$$\mu^{\ast}_s=\sum\limits_{j=1}^m\Big(\sum\limits_{p=1}^k M_{s,p}(\Delta_m)V_{p,j}(\Delta_m)
\Big)\lambda^{\ast}_j + \sum \limits_{p=1}^k
M_{s,p}(\Delta_m)g_{p}(f,\Delta_m) \eqno(3.7)$$ we determine the components
$\mu^{\ast}=(\mu^{\ast}_1,\mu^{\ast}_2,\ldots,\mu^{\ast}_k)\in
R^{nk}$ and construct the function
$$\mathcal{F}^{\ast}(t)=\sum\limits_{s=1}^k \varphi_s(t)\Big[\mu^{\ast}_s+\sum\limits_{r=1}^m
\int_{t_{r-1}}^{t_r}\psi_s(\tau)d\tau \lambda^{\ast}_r \Big]+f(t).
\eqno(3.8)$$ Recall that  $\lambda_r^{\ast}=x^{\ast}(t_{r-1}),$
where $x^{\ast}(t)$ is a solution to problem (2.1), ( 2.2). Hence, by solving system (3.6), we get the values of the desired solution at the left endpoints of the partition subintervals.

In order to determine the values of the function $x^{\ast}(t)$ at the remaining points of the subintervals $[t_{r-1},t_r)$, we solve the following Cauchy problems for the ordinary differential equation:
$$ \frac{dx}{dt} =
A(t)x + \mathcal{F}^{\ast}(t), \quad x(t_{r-1})=\lambda_r^{\ast},
\quad t\in[t_{r-1},t_r), \quad   r=\overline{1,m}.$$

Thus, the proposed algorithm contains seven interrelated parts.

If the fundamental matrices  $X_r(t),$ $r=\overline{1,m},$ are known, then the equalities   (2.16),(2.17), and (2.18) enable us to construct the system (3.6). Let $\lambda^{\ast}=(\lambda_1^{\ast},\lambda_2^{\ast},\ldots,\lambda_m^{\ast})\in
R^{nm}$ be a solution to (3.6). Then, using (3.7) and (3.8), we construct the function $\mathcal{F}^{\ast}(t)$ and determine a solution to problem (2.1), (2.2) by the equalities
$$x^{\ast}(t)=X_r(t)X_r^{-1}(t_{r-1})\lambda_r^{\ast}+X_r(t)\int_{t_{r-1}}^{t}
X_r^{-1}(\tau)\mathcal{F}^{\ast}(\tau)d\tau, \quad
t\in[t_{r-1},t_r), \quad r=\overline{1,m}, \eqno (3.9)$$
$$x^{\ast}(T)=X_m(T)X_m^{-1}(t_{m-1})\lambda_m^{\ast}+X_m(T)\int_{t_{m-1}}^{T}
X_m^{-1}(\tau)\mathcal{F}^{\ast}(\tau)d\tau. \eqno (3.10)$$ So, in this case the proposed algorithm provides the solution to the linear multipoint boundary value problem for integro-differential equations (2.1), (2.2) in the form  (3.9), (3.10).

As is known, it is not always possible to construct a fundamental matrix for a system of ordinary differential equations with variable coefficients. For this reason, we propose the following numerical implementation of the algorithm that is based on the fourth-order Runge-Kutta method and Simpson's rule.

I. Let us take a partition $\Delta_m:$
$t_0=0<t_1<\ldots<t_{m-1}<t_m=T$. We divide each $r$th subinterval $[t_{r-1},t_r],$   $r=\overline{1,m},$ into $m_r$ parts with the step size $h_r=(t_r-t_{r-1})/m_r.$

Suppose that on each subinterval $[t_{r-1},t_r]$  a variable $\widehat{t}$ takes on discrete values: $\widehat{t}=t_{r-1},$
$\widehat{t}=t_{r-1}+h_r,$ $\ldots,$
 $\widehat{t}=t_{r-1}+(m_r-1)h_r,$  $\widehat{t}=t_r$. Let $\{t_{r-1},t_r\}$ denote the set of such points.

II. Using the fourth-order Runge-Kutta method, we obtain numerical solutions to the Cauchy problems (3.1) and determine the values of the $(n\times n)$ matrix
$E_{\ast,r}^{h_r}(A(\cdot),\varphi_j(\cdot),\widehat{t})$ on the set $\{t_{r-1},t_r\}$, $r=\overline{1,m}$, $j=\overline{1,k}.$

III. Using the values of the $(n\times n)$ matrices $\psi_j(s)$ and
$E_{\ast,r}^{h_r}\Big(A(\cdot),\varphi_j(\cdot),\widehat{t}\Big)$ on
$\{t_{r-1},t_r\}$ and applying Simpson's rule, we determine the $(n\times n)$ matrices
$$\widehat{\psi}_{p,r}^{h_r}(\varphi_j)=\int _{t_{r-1}}^{t_r}\psi_p(\tau)E_{\ast,r}^{h_r}(A(\cdot),\varphi_j(\cdot),\tau)d\tau,
\quad p,j=\overline{1,k}, \quad  r=\overline{1,m}.$$

Summing up the matrices $\widehat{\varphi}_{p,r}^{h_r}(\psi_j)$ with respect to $r$, we get the $(n\times n)$ matrices
$G_{p,j}^{\widetilde{h}}(\Delta_m)=\sum\limits_{r=1}^m
\widehat{\varphi}_{p,r}^{h_r}(\psi_j)$, where
$\widetilde{h}=(h_1,h_2,\ldots,h_m)\in R^n$. We then construct the $(nk\times nk)$ matrix
$G^{\widetilde{h}}(\Delta_m)=(G_{p,j}^{\widetilde{h}}(\Delta_m)),$
$p,j=\overline{1,k}$.

Determine whether the matrix
$[I-G^{\widetilde{h}}(\Delta_m)]: R^{nk}\rightarrow R^{nk}$ is invertible. If so, we find
$[I-G^{\widetilde{h}}(\Delta_m)]^{-1}=(M_{p,j}^{\widetilde{h}}(\Delta_m)),$
$p,j=\overline{1,k}.$

In the case $[I-G^{\widetilde{h}}(\Delta_m)]$ is not invertible, we choose a new partition. In particular, as shown above, each subinterval can be divided in half.

IV. Solving the Cauchy problem (3.5), (3.6) by the fourth-order Runge-Kutta method, we get
 the values of the $(n\times n)$ matrix  $E_{\ast,r}(A(\cdot),A(\cdot),\widehat{t})$ and  $n$
vector $E_{\ast,r}(A(\cdot),f(\cdot),\widehat{t}) $ on
$\{t_{r-1},t_r\},$ $r=\overline{1,m}.$

V.  By using Simpson's rule on the grid  $\{t_{r-1},t_r\}$, we evaluate the definite integrals
$$\widehat{\psi}_{p,r}^{h_r}=\int_{t_{r-1}}^{t_r}\psi_p(\tau)d\tau,
\quad
\widehat{\psi}_{p,r}^{h_r}(A)=\int_{t_{r-1}}^{t_r}\psi_p(\tau)E_{\ast,r}^{h_r}
(A(\cdot),A(\cdot),\tau)
d\tau,$$
$$\widehat{\psi}_{p,r}^{h_r}(f)=\int_{t_{r-1}}^{t_r}\psi_p(\tau)E_{\ast,r}^{h_r}
(A(\cdot),f(\cdot),\tau)d\tau,  \quad r=\overline{1,N}, \quad
p=\overline{1,m}.$$

We the determine the $(n\times n)$ matrices
$V_{p,r}^{\widetilde{h}}(\Delta_m)$ and $n$ vectors  $g_{p}^{\widetilde{h}}(f,\Delta_m)$, $r=\overline{1,m}$, $p=\overline{1,k}$, by the equalities
$$V_{p,r}^{\widetilde{h}}(\Delta_m)=\widehat{\psi}_{p,r}^{h_r}(A)+
\sum\limits_{j=1}^{m}\sum\limits_{i=1}^{k}\widehat{\psi}_{p,j}^{h_j}(\varphi_i)
 \cdot\widehat{\psi}_{i,r}^{h_r},  \quad
g_{p}^{\widetilde{h}}(f,\Delta_m)=\sum\limits_{r=1}^{m}\widehat{\psi}_{p,r}^{h_r}(f).$$

VI. We construct the system of linear algebraic equations in parameters
$$Q_{\ast}^{\widetilde{h}}(\Delta_m) \lambda = -F_{\ast}^{\widetilde{h}}(\Delta_m),
 \qquad \lambda\in R^{nm}.  \eqno (3.11)$$
The elements of the matrix $Q_{\ast}^{\widetilde{h}}(\Delta_m)$ and the vector
$F_{\ast}^{\widetilde{h}}(\Delta_m)=(-d+ B_m F_{m}^{\widetilde{h}}(\Delta_m), F_{1}^{\widetilde{h}}(\Delta_m),
\ldots, F_{m-1}^{\widetilde{h}}(\Delta_m))$  are determined by the equalities
$$D_{r,i}^{\widetilde{h}}(\Delta_m)=\sum\limits_{j=1}^{k}E_{\ast,r}^{h_r}(A(\cdot),\varphi_j(\cdot),t_r)
\Big[\sum\limits_{p=1}^{k}M_{j,p}^{\widetilde{h}}(\Delta_m)V_{p,i}^{\widetilde{h}}(\Delta_m)
+\widehat{\psi}_{j,i}^{h_i}\Big],   \quad  i\neq r,   \quad
r, i=\overline{1,m},  $$
$$D_{r,r}^{\widetilde{h}}(\Delta_m)=\sum\limits_{j=1}^{k}E_{\ast,r}^{h_r}(A(\cdot),\varphi_j(\cdot),t_r)
\Big[\sum\limits_{p=1}^{k}M_{j,p}^{\widetilde{h}}(\Delta_m)V_{p,r}^{\widetilde{h}}(\Delta_m)
+\widehat{\psi}_{j,r}^{h_r}\Big]+E_{\ast,r}^{h_r}(A(\cdot),A(\cdot),t_r),
\quad r=\overline{1,m}, $$
$$F_{r}^{\widetilde{h}}(\Delta_m)=\sum\limits_{j=1}^{k}E_{\ast,r}^{h_r}(A(\cdot),\varphi_j(\cdot),t_r)
\sum\limits_{p=1}^{k}M_{j,p}^{\widetilde{h}}(\Delta_m)g_{p}^{\widetilde{h}}(\Delta_m)
+E_{\ast,r}^{h_r}(A(\cdot),f(\cdot),t_r), \quad
r=\overline{1,m}.$$

Using the constructed matrix
$(Q_{\ast}^{\widetilde{h}}(\Delta_m),$  we can establish the well-posedness of problem (2.1), (2.2). Suppose the matrix
$Q_{\ast}^{\widetilde{h}}(\Delta_m)$ is invertible and the estimate 
 $||Q_{\ast}(\Delta_m)-Q_{\ast}^{\widetilde{h}}(\Delta_m)||\leq
\varepsilon(\widetilde{h})$ holds.  If the inequality 
$||[Q_{\ast}^{\widetilde{h}}(\Delta_m)]^{-1}||\cdot
\varepsilon(\widetilde{h}) < 1 $ is true, then, by Theorem 4 [11, p.212], the matrix $Q_{\ast}(\Delta_m)$ is invertible. It follows then from Theorem 2.2 that problem (2.1),(2.2) is well-posed.

By solving (3.11) we determine $\lambda^{\widetilde{h}}\in
R^{nm}.$ As noted above, the elements
$\lambda^{\widetilde{h}}=(\lambda^{\widetilde{h}}_1,
\lambda^{\widetilde{h}}_2,\ldots,\lambda^{\widetilde{h}}_m)$ are the values of an approximate solution to problem (2.1), (2.2) at the left endpoints of the subintervals:
$x^{\widetilde{h}_r}(t_{r-1})=\lambda^{\widetilde{h}}_r,$
$r=\overline{1,m}.$

VII. In order to calculate the values of the approximate solution at the remaining points of the set $\{t_{r-1},t_r\}$,  we first find
$$\mu_i^{\widetilde{h}}=\sum\limits_{j=1}^{m}\Big(\sum\limits_{p=1}^{k}M_{i,p}^{\widetilde{h}}(\Delta_m)V_{p,j}
^{\widetilde{h}}(\Delta_m)\Big)\lambda_i^{h}+\sum\limits_{p=1}^{k}M_{i,p}^{\widetilde{h}}(\Delta_m)g_{p}
^{\widetilde{h}}(f,\Delta_m), \quad i=\overline{1,k}, $$ and then, using the fourth-order Runge-Kutta method, solve the Cauchy problems
$$ \frac{dx}{dt} =
A(t)x + \mathcal{F}^{\widetilde{h}}(t), \quad
x(t_{r-1})=\lambda_r^{\widetilde{h}}, \quad t\in[t_{r-1},t_r],
\quad r=\overline{1,m},$$ where
$$\mathcal{F}^{\widetilde{h}}(t)=\sum\limits_{i=1}^{k}\varphi_i(t)\Big(\mu_i^{\widetilde{h}}+
\sum\limits_{j=1}^{m}\widehat{\psi}
_{i,j}^{{h}_j}\lambda_j^{h}\Big)+f(t).$$

Thus the algorithm allows us to find a numerical solution to problem (2.1), (2.2).

\vskip0.5cm

{\bf 4.} {\bf A multipoint problem for an integro-differential equation with degenerate kernel}

Let us now turn to the original multipoint problem (1.1), (1.2). To solve the problem, we will approximate the kernel of the integral summand by a degenerate kernel [5, 6, 8, 10].

By the Weierstrass polynomial approximation theorem, for any $\varepsilon>0$ there exist a number  $k= k(\varepsilon)$ and continuous on $[0,T]$ matrices $\varphi_j(t),$ $\psi_j(\tau),$
$j =\overline{1,k}$, such that the following inequality holds
$$\max\limits_{t\in[0,T]}\int\limits_{0}^T \|K(t,\tau)-\sum\limits_{j=1}^{k}
\varphi_j(t)\psi_j(\tau)\|d\tau <\varepsilon . \eqno(4.1)$$

The set of matrices $\{\varphi_j(t), \psi_j(\tau), j=\overline{1,m}\},$ satisfying (4.1), we will call the $\varepsilon$-approximating set for $K(t,\tau)$. The multipoint problem with degenerate kernel (2.1), (2.2), corresponding to (1.1),(1.2), we will call the $\varepsilon$-approximating problem for problem (1.1), (1.2).

Assuming the $\varepsilon$-approximating multipoint problem (2.1), (2.2) to be well-posed with constant $C_k$, we find the solution to problem (1.1), (1.2) according to the following algorithm.

{\bf Step 0.} By solving problem (2.1), (2.2) we get a function $x^{(0)}(t)$, which we take as an initial approximation to the solution to problem (1.1), (1.2).

{\bf Step 1.}  Using $x^{(0)}(t)$ and solving the $\varepsilon$-approximating problem
$$ \frac{dx}{dt} =
A(t)x + \sum\limits_{j=1}^k\varphi_j(t)\int ^T_0
\psi_j(\tau)x(\tau)d\tau +f(t)+$$ $$+\int ^T_0 \Bigl[K(t,\tau)-\sum\limits_{j=1}^k\varphi_j(t)\psi_j(\tau)\Bigr]x^{(0)}(\tau)d\tau,
\qquad t\in (0,T), \eqno (4.2)$$
$$ \sum \limits ^m_{i=0}B_i x(t_i) = d,  \qquad  d \in R^n, \eqno (4.3)$$
we get the function $x^{(1)}(t)$.

Continue the algorithm, in the $i$th step $(i=1,\ldots)$ we solve the problem
$$ \frac{dx}{dt} =
A(t)x + \sum\limits_{j=1}^k\varphi_j(t)\int ^T_0
\psi_j(\tau)x(\tau)d\tau +f(t)+$$ $$+\int
^T_0\Bigl[K(t,\tau)-\sum\limits_{j=1}^k\varphi_j(t)\psi_j(\tau)\Bigr]x^{(i-1)}(\tau)d\tau,
\qquad t\in (0,T), \eqno (4.4)$$
$$ \sum \limits ^m_{i=0}B_i x(t_i) = d, \eqno (4.5)$$
and get the function $x^{(i)}(t)$.

The well-posedness of the approximating problem ensures the feasibility of the algorithm and allows us to construct the sequence $(x^{(i)}(t)),$ $i=0,1,\ldots.$

The following assertion provides conditions for the convergence of the algorithm to the unique solution of multipoint problem (1.1), (1.2) and the estimates for the difference between the exact and approximate solutions to the problem.

{\bf Theorem 4.1} Let the $\varepsilon$-approximating problem (2.1), (2.2) be well-posed with constant $C_k$. Suppose that the following inequality holds:
$$ q_m^{\varepsilon}=K_m\cdot \varepsilon<1.  \eqno (4.6) $$
Then the algorithm converges to $x^{\ast}(t)$ and the estimate
$$\|x^{\ast}-x^{(i)}\|_1\leq\frac{1}{1-q_k^{\varepsilon}}(q_k^{\varepsilon})
^{i}\cdot C_k\max(\|f\|_1,\|d\|)  \eqno (4.7)$$ is valid, where
$x^{\ast}(t)$ and $x^{(i)}(t)$ are the unique solutions to problems (1.1), (1.2) and (4.4), (4.5), respectively.

{\bf Proof.} By assumption, there exists a unique solution to problem (2.1),(2.2) and it satisfies the inequality
$$\|x^{(0)}\|_1\leq C_k \max(\|f\|_1,\|d\|).$$
By solving problem (4.2), (4.3), we get $x^{(1)}(t).$  The difference $\Delta x^{(1)}(t)= x^{(1)}(t)-x^{(0)}(t)$ satisfies the inequality
$$\|\Delta x^{(1)}\|_1\leq C_k \cdot
\max\limits_{t\in[0,T]}\int_{0}^{T}
\|K(t,s)-\sum\limits_{j=1}^k\varphi_j(t)\psi_j(\tau)\|d\tau
\|x^{(0)}\|_1\leq C_k \cdot\varepsilon\cdot
C_k\cdot\max(\|f\|_1,\|d\|).  \eqno (4.8)$$ 

In the same way, by solving problem (4.4), (4.5), we get $x^{(i)}(t)$, and for the difference
$\Delta x^{(i)}(t)= x^{(i)}(t)-x^{(i-1)}(t)$ we have
$$\|\Delta x^{(i)}\|_1\leq C_k\cdot
\varepsilon\cdot \|\Delta x^{(i-1)}\|_1=q_k ^{\varepsilon}\|\Delta
x^{(i-1)}\|_1, \quad  i=2,3,\ldots. \eqno (4.9)$$ 

The convergence of the sequence $(x^{(i)}(t)),$ $i=0,1,\ldots$ to the solution $x^{\ast}(t)$ of problem (1.1), (1.2), as well as the uniqueness of this solution, follow from the inequalities (4.6) and (4.9). The estimate is derived from (4.8) and (4.9).

Theorem 4.1 is proved.

The conditions of Theorem 2.1 ensure the existence of a unique solution to problem (2.1), (2.2.) and the validity of the estimate (2.23). The number $\mathcal{N}(k,\Delta_m)$ in (2.23), as mentioned above, does not depend on  $f(t)$ and $d$. We therefore can treat this number as a constant of well-posedness of problem (2.1), (2.2). Hence, by Theorems 2.1 and 4.1, the following statement holds true.

{\bf  Theorem 4.2} Suppose that

(a) the set of the matrices $\{\varphi_j(t), \psi_j(\tau),
j=\overline{1,k}\}$ is an $\varepsilon$-approximating set for $K(t,\tau)$;

(b) $\Delta_m\in \sigma(k,[0,T]);$

(c) the matrix $Q_{\ast}(\Delta_m): R^{nm}\rightarrow R^{nm}$ in (2.22) is invertible;

(d) the inequality 
$\delta_m^{\varepsilon}=\mathcal{N}(k,\Delta_m)\cdot
\varepsilon<1$ holds.

Then problem (1.1), (1.2) is well-posed with constant
$\displaystyle{C=\frac{1}{1-\delta_k^{\varepsilon}}\cdot\mathcal{N}(k,\Delta_m)}.$

The conditions of Theorem 4.1 are not only necessary but also sufficient for the well-posedness of problem (1.1), (1.2).

{\bf Theorem 4.3} Problem (1.1), (1.2) is well-posed if and only if there exists $\varepsilon$-approximating multipoint problem (2.1),(2.2), that is well-posed with constant $C_k$, and the inequality (4.6) holds true.

{\bf Proof.} The sufficiency of the conditions of Theorem for the well-posedness of problem (1.1), (1.2) follow from Theorem 4.2.

Let us prove the necessity. Assume that problem (1.1), (1.2) is well-posed with a constant $C$. Take $\varepsilon>0$ satisfying the inequality $\varepsilon \cdot C<1/2.$ For chosen $\varepsilon$ take $k \in \mathbb{N}$ and continuous on $[0,T]$ matrices $\varphi_j(t),$ $\psi_j(\tau),$
$j = \overline{1,k},$ satisfying inequality (4.1). Let us show that multipoint problem (2.1), (2.2) with these matrices is well-posed and the constant $C_k$ of well-posedness satisfies inequality (4.4). To this end, we use the following algorithm.

{\bf Step 0.} By solving problem (1.1), (1.2), we get the function $x^{(0)}(t).$

{\bf Step i.}  Assuming $x^{(i-1)}(t),$  $i=1,2,\ldots,$ to be known, we solve the problem
$$ \frac{dx}{dt} =
A(t)x + \int  ^T_0 K(t,\tau)x(\tau)d\tau +f(t)+ \int
^T_0\Bigl[\sum\limits_{j=1}^k \varphi_j(t)\psi_j(\tau)-K(t,\tau)\Bigr]x^{(i-1)}(\tau)d\tau,
\qquad t\in [0,T], $$
$$ \sum \limits ^m_{i=0}B_i x(t_i) = d,$$ and get the function  $x^{(i)}(t).$

It is easy to check that the algorithm converges to $x^{\ast}(t)$ and the estimate
$$\|x^{\ast}\|_1\leq\frac{C}{1- C\cdot\varepsilon}\cdot \max(\|f\|_1,\|d\|),   \eqno (4.10)$$
holds, were $x^{\ast}(t)$ is the unique solution to problem (2.1), (2.2).

Since, by assumption, $C\cdot\varepsilon<1/2,$ the well-posedness of the $\varepsilon$-approximating problem (2.1),(2.2) with constant $C_k=2C$ follows from (4.10). Taking into account the choice of $\varepsilon>0$, we get $q_k^{\varepsilon} = C_k\cdot \varepsilon < 1.$

Theorem 4.3 is proved.

\vskip0.5cm

 {\bf  References}

\vskip0.3cm {\footnotesize \noindent {[1]} {V. M. Abdullaev, K. R.
Aida-zade},  {\it Numerical method of solution to loaded nonlocal
boundary value problems for ordinary differential equations},
Comput. Math. Math. Phys., {\bf 54}(2014), pp. 1096-1109.\\
{[2]} {K. R. Aida-zade,   V. M. Abdullaev},  {\it On the numerical
solution of loaded systems of ordinary differential equations with
nonseparated multipoint and integral conditions},  Numer. Anal. Appl., {\bf 7}(2014), pp. 1-14.\\
{[3]} {A. T. Assanova, E. A. Bakirova, Zh. M. Kadirbayeva}
{\it Numerical solution to a control problem for integro-differential equations},
Comput. Math. and Math. Phys., {\bf 60} (2020), pp. 203--221. \\
{[4]} {A. T. Assanova, A. E. Imanchiyev, Z. M. Kadirbayeva}
{\it Numerical solution of systems of loaded ordinary differential
equations with multipoint conditions},
Comput. Math. and Math. Phys., {\bf 58} (2018), pp. 508-516. \\
{[5]} {K.I. Babenko},  {\it Fundamentals of Numerical Analysis},
 Nauka, 1986. [in Russian] \\
{[6]}  {N.S. Bakhvalov}, {\it Numerical Methods},
 Fizmatgiz, 1973.   [in Russian] \\
{[7]} {A. A. Boichuk, A. M. Samoilenko}, {\it Generalized inverse
operators and Fredholm boun\-dary-value problems,} VSP,
Utrecht,  2004. \\
{[8]}  {H. Brunner}, {\it Collocation methods for Volterra
integral
and related functional equations,} Cambridge University Press, 2004.   \\
{[9]} {Ya. V. Bykov,} On Some Problems in
the Theory of Integro-Differential Equations. Kirgiz. Gos. Univ.,
Frunze, 1957. [in Russian]\\
{[10]} H. Cohen, {\it Numerical approximation methods,}  Springer, 2011.  \\
{[11]}  L.M. Delves,  J.L. Mohamed,  {\it Computational methods
for integral equations,} Cambridge University Press, 1985. \\
{[12]}  D.S. Dzhumabayev,  {\it Criteria for the unique
solvability of a linear boundary-value problem for an ordinary
differential
equation}, U.S.S.R. Comput. Maths. Math. Phys., {\bf 29}(1989), pp. 34-46.\\
{[13]}  D.S. Dzhumabaev, {\it A method for solving the linear
boundary value problem for an integro-differential equation},
Comput. Math. Math. Phys., {\bf 50}(2010), pp. 1150-1161.\\
{[14]}  D.S. Dzhumabaev, {\it Necessary and sufficient conditions
for the solvability of linear boundary-value problems for the
Fredholm
integro-differential equation},  Ukr. Math. J., {\bf 66}(2015), pp. 1200-1219.\\
{[15]}  D.S. Dzhumabaev, {\it An algorithm for solving the linear
boundary value problem for an integro-differential equation},
Comput. Math. Math. Phys., {\bf 53}(2013), pp.  736-758.\\
{[16]}  D.S. Dzhumabaev, {\it On one approach to solve the linear
boundary value problems for Fredholm integro-differential
equations},  J. Comput. Appl. Math., {\bf 294}(2016), pp. 342-357. \\
{[17]} D. S. Dzhumabaev, {\it Computational methods of solving the boundary value problems for the loaded differential and Fredholm integro-differential equations}, Math. Meth. Appl. Sci., {\bf 41} (2018), pp. 1439-1462.\\
{[18]} D. S. Dzhumabaev, {\it New general solutions to linear Fredholm integro-differential equations and their applications on solving the boundary value problems}, J. Comput. Appl. Math., {\bf 327}(2018), pp. 79-108. \\
{[19]} D. S. Dzhumabaev, E. A. Bakirova, {\it Criteria for the
unique solvability of a linear two-point  boundary value problem
for systems of integro-differential equations},  Differ.
Equ., {\bf 49}(2013), pp.914-937. \\
{[20]} D. S. Dzhumabaev, A. E. Imanchiev, {\it Well-posed solvability of linear multipoint boundary value problem}, Mathematical Journal., {\bf 5}:1 (2005), pp. 30-38. [in Russian] \\
{[21]} A. E. Imanchiev, {\it Necessary and sufficient conditions of the unique solvability of linear multipoint boundary value problem}, News of MES RK, NAS RK. Phys.-Mathem. Series.,   {\bf 3} (2002), pp. 79-84. [in Russian] \\
{[22]}  K. Maleknejad,  M. Attary, {\it An efficient numerical
approximation for the linear Fredholm integro-differential
equations based on Cattani's method},  Commun. Nonlinear Sci. Numer. Simulat., {\bf 16}(2011), pp. 2672-2679.\\
{[23]}  A. Molabahrami, {\it Direct computation method for solving
a general nonlinear Fredholm integro-differential equation under
the mixed conditions: Degenerate and non-degenerate
kernels},  J. Comput. Appl. Math., {\bf 282}(2015), pp. 34-43.\\
{[24]} \c{S}.  Yuzbasi, {\it Numerical solutions of system of
linear Fredholm-Volterra integro-differential equations by the
Bessel collocation method and error estimation}, Appl. Math. Comput., {\bf 250}(2015), pp. 320-338.\\
{[25]} A.M. Wazwaz, {\it  Linear and Nonlinear Integral
Equations: Methods and Applications,} Higher Education Press,
Beijing and Springer-Verlag, 2011.}

\end{document}